\documentclass[a4paper]{article}

\usepackage{url}
\usepackage{graphicx}
\usepackage{subcaption}
\usepackage{amssymb}
\usepackage{amsthm}
\usepackage{amsmath}
\usepackage{geometry}

\usepackage{algorithm}
\usepackage{pgfplots}

\usepackage{algpseudocode}
\usepackage{tikz}
\usepackage{verbatim}

\usepackage{hyperref}
\usepackage{authblk}

\usetikzlibrary{shapes,positioning,arrows,calc,shapes.geometric}
\usetikzlibrary{positioning}
\usetikzlibrary{calc,patterns,positioning,decorations.pathmorphing,arrows,decorations.markings} 
\usepgfplotslibrary{groupplots}
\pgfplotsset{compat=newest}

\usepackage{float}
\usepackage{mathrsfs}

\newcommand{\iu}{\mathrm{i}} 
\newcommand{\cF}{\mathcal{F}}
\newcommand{\FT}[1]{\cF\left\{ #1 \right\}}
\newcommand{\IFT}[1]{\cF^{-1}\left\{ #1 \right\}}

 \usepackage[figcolor=white]{todonotes}

\begin{document}

\markboth{D. Walsken, P. Petrov and M. Ehrhardt}{A Spectral Split-Step Padé Method for Guided Wave Propagation}
\title{A Spectral Split-Step Padé Method for Guided Wave Propagation}

\author[1]{Daniel Walsken\footnote{\tt walsken@math.uni-wuppertal.de}}
\author[1]{Matthias Ehrhardt\footnote{\tt ehrhardt@math.uni-wuppertal.de}}
\author[2]{Pavel Petrov\footnote{\tt pavel.petrov@impa.br}}

\affil[1]{University of Wuppertal, Chair of Applied and Computational Mathematics, Gauß\-straße 20, 42119 Wuppertal, Germany}
\affil[2]{Instituto de Matemática Pura e Aplicada, Estrada Dona Castorina~110, Rio de Janeiro, CEP 22460-320, Brazil}

\maketitle

\begin{abstract}
    In this study, a Fourier-based, split-step, Padé (SSP) method for sol\-ving the parabolic wave equation with applications in guided wave propagation in ocean acoustics is presented.
    Traditional SSP implementations rely on finite-difference discretizations of the depth-dependent differential operator. 
    This approach limits accuracy in coarse discretizations as well as computational efficiency in dense discretizations, since it is does not significantly benefit from parallelization. In contrast, our proposed method replaces finite differences with a spectral representation using the discrete sine transform (DST). 
    This enables an exact treatment of the vertical operator under homogeneous boundary conditions. 
    For non-constant sound speed profiles, we use a Neumann series expansion to treat inhomogeneities as perturbations. 
    Numerical experiments demonstrate the method’s accuracy in range-independent and range-dependent scenarios, including propagation in deep ocean with Munk profile and in the presence of a parameterized synoptic eddy.
    Compared to finite-difference SSP methods, the Fourier-based approach achieves higher accuracy with fewer depth discretization points and avoids the resolution bottleneck associated with sharp field features, making it 
    
    well-suited for large-scale, high-frequency wave propagation problems in ocean environments.
\end{abstract}

\section{Introduction}\label{S1}
The parabolic equation method originates from the work of Leontovich and Fock \cite{leontovich1946,vlasov1995parabolic} where it was proposed as a practical tool for radiowave propagation over the Earth's surface.

This formulation reduced the computational burden by transforming the full-wave elliptic equation into an evolutionary equation that can be efficiently solved by a one-way marching numerical scheme (with back scattering neglected). In geophysics, Claerbout applied similar principles to seismic waves simulation \cite{claerbout1976fundamentals}, while Tappert later on 
introduced the parabolic equation in ocean acoustics \cite{Tappert1977}.

Paraxial (narrow-angle) parabolic equations from early papers were not capable of handling the waves propagating at large grazing angles to the waveguide axis. For instance, in ocean acoustics this manifests in errors when modelling sound reflection from the seabottom. This issue was resolved by introducing the so-called wide-angle parabolic equations \cite{popov1977,claerbout1976fundamentals} (WAPEs) which are currently widely used in optics \cite{lu2006some}, radio waves theory \cite{lytaev2019} and atmospheric and underwater acoustics \cite{collins93,abawi1997coupled,ostashev2019extra}.

Subsequent developments led to the creation of mode parabolic equations describing the evolution of amplitudes of normal modes in a three-dimensional waveguide. The adiabatic mode parabolic equation (PE) \cite{collins1993adiabatic,petrov2012,petrov2020pseudodifferential} assumes weak mode coupling, whereas the coupled mode PEs  \cite{abawi1997coupled,trofimov2015mode,petrov2024generalization,he2024novel} take mode interactions into account.

The WAPE solution technique known as \emph{split-step Pad\'e} (SSP) method, independently proposed by Collins~\cite{collins93} and Avilov~\cite{avilov1985} has become an essential 
step in the development of the parabolic equation theory.

In this approach, the Pad\'e series is used to approximate the exponential of the square-root of a differential operator  (rather than the square root itself) in the transverse direction to the waveguide axis. A finite-difference discretization of this differential operator is commonly used to compute this Pad\'e approximation numerically. Such discretization imposes certain restrictions on the meshsize and introduces a truncation error for the derivatives in this direction.

Modern spectral methods overcome this limitation by precisely representing the vertical differential operator in spectral space, allowing for accurate simulations on coarser grids.
The first such method, the \textit{split-step Fourier} (SSF) method, introduced independently by both Hardin~\cite{hardin1973applications} and Tappert~\cite{tappert1974numerical}, is suitable for narrow-angle parabolic equations, that is, for equations obtained from first-order Taylor expansion. 
The aim of this study is to combine the high accuracy at high grazing angles of the SSP method with the exact representation of the differential operator in the SSF method.

Compared to solving the full wave equation numerically (e.g., \cite{Specfem3d}), the PE approach requires far less computational effort 
because direct discretizations must finely resolve each wavelength. Helmholtz equation solvers based on separation of variables and finite differences 
(e.g., \cite{porter1992kraken, orca2002}) also incur a high computational cost due to coupled eigenvalue problems. 
In both cases, convergence depends on the spatial discretization because the depth operator is approximated by finite differences. 
In constant sound-speed (isovelocity) environments, however, Fourier transforms apply the operator exactly in spectral space, avoiding this restriction entirely. 
This improves both accuracy and 
also efficiency by making the propagation step inherently parallelizable.

The remainder of this paper is organized as follows:  
In Section~\ref{S2}, we discuss the split step Pad\'e (SSP) method in its original form (hereafter abbreviation SSP is reserved for its traditional implementation based the on finite-difference discretization). In Section~\ref{S3}, we introduce the Fourier Transform as a tool to solve the parabolic wave equation on isovelocity domains and extend the method to also solve equations with variable sound speed. Section~\ref{S4} presents numerical results demonstrating the method's accuracy and efficiency.
Finally, Section~\ref{SC} concludes with a discussion of the findings.

\section{The Split-Step Pad\'e Method}\label{S2}
The parabolic wave equation is commonly derived from the two-dimensional Helmholtz equation
\begin{equation}\label{HE}
    \partial_r^2 u + \partial_z^2u + k^2u = 0\,,
\end{equation}
where $k = k(r,z) = \tfrac{\omega}{c(r,z)}$ is the wavenumber and $u$ is the complex acoustic field. Note that Eq.~\eqref{HE} can describe both a line source in plane geometry and a point source in the cylindrical geometry \cite{jensen2011computational} (with the spreading factor $1/\sqrt{r}$ removed).

By formally factoring the elliptic differential operator into
\begin{equation}
    \Bigl(\partial_r + \iu \sqrt{\partial_z^2 + k^2}\Bigr)\,\Bigl(\partial_r - \iu \sqrt{\partial_z^2 + k^2}\Bigr)\,u = 0
\end{equation}
one can obtain separate equations for left- and right-propagating waves, respectively. Neglecting back scattering effects (which is a common approximation for guided wave propagation) we restrict our attention to the propagation in the positive $r$ direction
\begin{equation}\label{PDPE}
    \partial_r u = \iu \sqrt{\partial_z^2 + k^2} \,u\,.
\end{equation}
The so-called pseudo-differential parabolic wave equation can then be formally solved by 
integration over the interval $[r_0, r_0+h]$, resulting in
\begin{equation}
     u(r+h,z) = \exp\Bigl(\iu h \sqrt{\partial^2_z +k^2}\Bigr)\,u(r,z) = \hat{P}u(r,z)\,,
\end{equation}
where the operator $\hat{P}$ is called propagator. 
Note that existence, uniqueness and well-posedness results for Eq.~\eqref{PDPE} were recently established \cite{ehrhardt2025} (a rigorous mathematical definition of the square root operator in the framework of a suitable operator calculus was also given in this study).

The idea of the \textit{split-step Pad\'e} (SSP) method \cite{collins93,petrov2020pseudodifferential} consists in using a Pad\'e approximation of the order $[p/p]$ for the propagator $\hat{P}$ as follows
\begin{equation}
    \exp\Bigl( \iu h \sqrt{\partial_z^2 + k^2}\Bigr)\,u 
    \approx e^{\iu k_0 h} \biggl( d_0 + \sum_{j=1}^{p}\frac{d_j}{1 + b_j \hat{X}} \biggr)\,u
    = e^{\iu k_0 h} \biggl( d_0u + \sum_{j=1}^{p}d_j w_j \biggr)\,
\end{equation}
with $\hat{X} = \tfrac{k^2 - k_0^2 + \partial_z^2}{k_0^2} = \frac{\delta k^2 + \partial_z^2}{k_0^2}$.
The weights $b_j$, $d_j$, $j=1,\ldots,p$ are the Pad\'e coefficients computed according to \cite[p.~27--38]{baker1975essentials}, while the functions $w_j$ are defined by the relation
\begin{equation} \label{eq:1dhelmholtz}
    \left(1 + b_j\hat{X}\right)\, w_j = u(r,z)\,.
\end{equation}
Because of the dependence of $\delta k^2$ on $z$, SSP technique usually relies on finite difference methods (FDMs) to approximate the operator $\partial_z^2$ acting on $u(r,z)$.
Note that if the medium wavenumber is constant (i.e., Eq.~\eqref{HE} has constant coefficients), then the differential operator can be represented by multiplication in Fourier domain analytically (which allows to resolve each Eq.~\eqref{eq:1dhelmholtz} efficiently using FFT).

\section{The Spectral Split-Step Pad\'e Method}\label{S3}

Denoting $\tilde{b}_j:=b_j/k_0^2$, we can rewrite Eq.~\eqref{eq:1dhelmholtz} on a finite interval $z\in[0,H]$ in the form
\begin{subequations}
    \begin{align}
        \Bigl( 1 +\tilde{b}_j \partial_z^2 +  \tilde{b}_j\delta k^2(z) \Bigr)\, w_j(z) &= u(r,z) \\
        \Bigl( \hat{L}_j + \tilde{b}_j\delta k^2(z) \Bigr) \,w_j(z) &= u(r,z) \,. \label{eq:1dHelmholtz}
    \end{align}
\end{subequations}
The part dependent on $z$, $\delta k^2$ can then later be regarded as a perturbation.
Consider now the isovelocity case
\begin{equation}
    \hat{L}_j w_j(z) = \Bigl(1 + \tilde b_j \delta k^2(z)\Bigr) w_j(z)= v(z)\,.
\end{equation}
Discretizing above equation along the $z$ axis with the equidistant grid $z_n=hn; \; h=\tfrac{H}{N}; \; n=0,\ldots,N$ and applying the \textit{discrete Fourier transform} (DFT) gives
\begin{equation}
    \left( 1 + \left( \ell\frac{2\iu \pi}{H}\right)^2\right) (\bar w_j)_\ell = \bar v_\ell \,
\end{equation}
for all $\ell$.
The lower index $\ell$ denotes the $\ell$th element of the discrete, Fourier transformed functions $\bar w_j = \FT{w_j}$ and $\bar v = \FT{v}$.
It is now straightforward to invert the operator analytically by the following formulae
\begin{subequations}
\begin{align}
    (\bar w_j)_\ell &= \left(1 + \left( \ell\frac{2\iu \pi}{H}\right)^2\right)^{-1} \bar v_\ell \\
    w_j &= \hat L_j^{-1} v = \IFT{\left(1 + \left( \ell\frac{2\iu \pi}{H}\right)^2\right)^{-1} \FT{v} }  \,, \label{eq:inverse_isovel}
    \end{align}
\end{subequations}
where $\IFT{\circ}$ denotes the inverse discrete Fourier transform.
Now considering again the full operator and regarding $\tilde b_j \delta k^2(z)$ as a small perturbation, i.e. $L_j w_j> \delta k^2 w_j$, we can rewrite Eq.~\eqref{eq:1dHelmholtz} as
\begin{equation}
    w_j(z)=\left( 1 - \left(-\tilde{b}_j \hat{L}_j^{-1} \delta k^2(z) \right) \right)^{-1}\hat{L}_j^{-1} u(r,z) \,.
\end{equation}
This solution operator can now be expressed as a Neumann series
\begin{subequations}
\label{eq:neumannseries}
\begin{align}
    w_j(z) &= \sum_{m=0}^{\infty}  \Bigl(- \tilde{b}_j\hat{L}_j^{-1}\delta k^2(z) \Bigr)^m \hat{L}_j^{-1} u(r,z) \\
    &\approx \sum_{m=0}^{M}  \Bigl(- \tilde{b}_j\hat{L}_j^{-1}\delta k^2(z) \Bigr)^m \hat{L}_j^{-1} u(r,z) \,, \label{eq:Neumann}
\end{align}
\end{subequations}
with $M\in\mathbb{N}$ being the cutoff index of the series. Here, the inverse isovelocity operator $\hat{L}_j^{-1}$ can be applied via Eq.~\eqref{eq:inverse_isovel} and the operator $\delta k^2$ can be applied via multiplication.
Note that $\hat{L}_j^{-1}$ and $\delta k^2$ are noncommutative.
As such, the order has to be taken into account, thus to compute the $m$th term of the series, a total of $m+1$ Fourier- and inverse Fourier Transforms have to be computed.

The computations described above can be summarized in the form of an algorithm for solving the parabolic wave equation. 
Let for this algorithm $N\in\mathbb{N}$ denote the number of steps to be taken, $b_j$, $d_j$, $j=1,\ldots,p$ denote the Pad\'e coefficients and $M+1$ be the number of series terms taken into account.

\begin{algorithm}[ht!]
    \caption{Spectral Split Step Pad\'e (SSSP)}
    \begin{algorithmic}
    \Require $u_0 = \{ u(0, n\Delta z )\colon n=0,\ldots,N  \}$; $h$
    \State {compute $d_j$ and $b_j$ for $j=0,\ldots,p-1$ according to \cite[p.~27--38]{baker1975essentials}}
        {\For{$n=1:N$}
        {\For{$j=0:p$}
        \State  $\tilde{b}_j \gets \tfrac{b_j}{k_0^2}$
        \State $w_j \gets \sum_{m=0}^{M}(-\hat{L}_j^{-1} \tilde{b}_j \delta k^2)^m \hat{L}_j^{-1} u_{n-1}$ 
        \EndFor}
        \State $u_n \gets e^{\iu h k_0}\sum_{j=0}^{p} d_j w_j$
        \EndFor}
\end{algorithmic}
\end{algorithm}

\subsection{Boundary Conditions}
Up to this point, boundary conditions have not been mentioned yet.
A major drawback of spectral methods in general and this method specifically is that it requires both the top and the bottom of the domain to satisfy the same boundary condition.
Take the same depth domain as introduced before, $z \in [0,H]$ and let the range domain be $r \in [0, R]$.
In the simplest case, as derived above, one assumes periodic boundary conditions
\begin{equation}
    u(r, 0) = u(r,H) \,.
\end{equation}
Those are imposed implicitly by the discrete Fourier transform (DFT).
In the case of more realistic homogeneous Dirichlet boundary conditions
\begin{equation}
    u(r,0) = u(r, H) = 0 \,,
\end{equation}
the DFT has to be replaced by a \textit{discrete sine transform} (DST), which satisfies the boundary, as well as for homogeneous Neumann boundary conditions
\begin{equation}
    \partial_z u(r, 0) = \partial_z u(r, H) = 0
\end{equation}
a discrete cosine transform can be employed.
The following numerical tests employ a DST to satisfy homogeneous Dirichlet BCs, 
but due to the nature of the chosen example problem could also use a DCT with Neumann BCs for similar results.

Of particular practical interest are transparent or absorbing BCs, as they allow for the truncation of the computational domain.
This allows computing fields on subsets of unbounded domains numerically.
While transparent BCs have not yet been developed for the presented method, a more generally applicable method, \textit{perfectly matched layers} (PML) will be considered in future work.
Another option would be to use complex absorbing potentials, which are not numerically optimal, but can be used as a simple replacement.
As the name suggests, this introduces an imaginary component in the wavenumber, which in turn acts as an absorbing layer around the boundaries.

\subsection{Initial Conditions}

The initial condition $u(0,z)=u_0(z)$ (often called the starter \cite{jensen2011computational})
for parabolic equations can be generated in a number of ways.

In most cases in acoustics, optics and radio physics, the source of the wave can be approximated as a point source.
To start the propagation, one can take the normal modes solution of this point source at $r=r_0$, $z=z_s$
\begin{equation}
    u(r_0,z) = \sqrt{2\pi} \sum_{m=1}^M \frac{\Psi_m(z_s) \Psi_m(z)}{\sqrt{\kappa_{m}}}\,.
\end{equation}
Those mode functions $\Psi_m$ as well as the wavenumbers $\kappa_m$ are obtained by solving the Sturm-Liouville problem
\begin{equation}
    \partial_z^2 \Psi_m + \left[ k^2(z) + \kappa^2 \right] \Psi_m = 0 \,
\end{equation}
on a discrete grid along the $z$ axis with appropriate boundary conditions.
In an acoustic deep water setting, homogeneous Dirichlet boundary conditions can be used, as the interaction between the pressure field and the boundaries is minimal.

A less expensive numerical starter is the self-starter developed by Collins~\cite{collins1992self}, 
which, while not used in this study, can in principle be applied with the Spectral Split Step Pad\'e (SSSP) approach with only minor modifications. 

Analytical starters have also been developed to accurately match the far field behavior of a point source.
Most notably, starters in the form of a Gaussian function are used. For wide angle PEs specifically, Greene's source~\cite{greene1984rational}
\begin{equation}
    u(r_0, z) = \sqrt{k_0}\left[ 1.4467 - 0.4201 k_0^2 (z - z_s)^2 \right] e^{- \frac{k^2 (z - z_s)}{3.0512}} \,, \label{eq:Greene}
\end{equation}
has been shown to have good spectral properties and is much simpler to compute than a modal starter.

A ray-based starter has been developed recently~\cite{petrov2020wide} and will be considered upon further development of SSSP.

\section{Numerical Results}\label{S4}

Since the convergence rate of the series \eqref{eq:neumannseries} depends on the relative magnitude of the wave-number perturbation $\delta k^2$, the first scenario focuses on the propagation of a time-harmonic acoustic signal in a range-independent deep-water waveguide (with a SOFAR channel) with perfectly reflecting surface and bottom. The second numerical experiment demonstrates the range-dependent capabilities of the methods by simulating sound propagation through a parametrized synoptic eddy. For the final experiment, we compare the performance of the classical SSP method with the newly developed SSSP method for very coarse depth discretizations. In this case, the error from approximating the differential operator by finite differences in the standard SSP is so large that the SSSP outperforms the classical method.
The software used can be found at~\cite{sssp}.

\subsection{Range-Independent Propagation}
In the first numerical experiment, we study the propagation of acoustic waves in a range-independent, deep-water environment using the Munk sound speed profile.
The objective of this subsection is to determine the number of terms in the Neumann series necessary to accurately reproduce the field using the SSSP.
The Munk profile 
\begin{equation}
    c_{\rm Munk}(z) = 1500 \,\bigl(1 + 0.00737(\eta + e^{-\eta} -1)\bigr)
    ~{\rm m}{\rm s}^{-1}; \quad \eta = 2\,\frac{z - 1300~{\rm m}}{1300~{\rm m}} \label{eq:Munk}
\end{equation}
is a widely used model for the dependence of sound speed on depth in the deep sea. 
We set the total depth to 4~km and use a normal mode-based starter \cite{jensen2011computational} representing a point source deployed at a depth of $z_s=1100~\mathrm{m}$,
so that the acoustic field at $r=0$ is computed as a linear combination of normal mode eigenfunctions. 
We compare our solution to the solution obtained by the classical finite-difference SSP method (i.e., SSP based on a finite-difference approximation of derivatives with respect to depth $z$).

\begin{figure}[ht!]
    \centering
    \begin{subfigure}{0.48\textwidth}
        \centering
        \includegraphics[width=\linewidth]{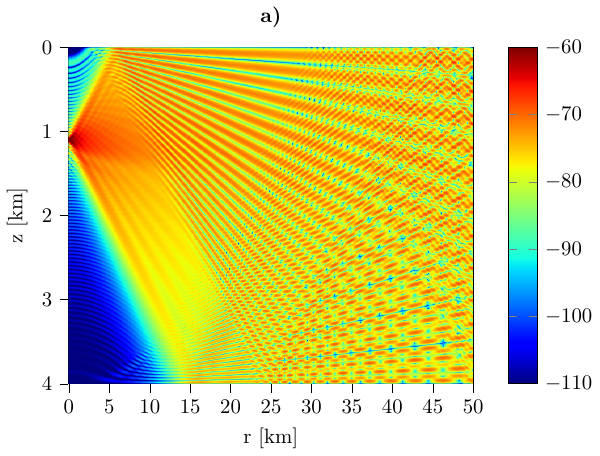}
    \end{subfigure}
    \hfill
    \begin{subfigure}{0.48\textwidth}
        \centering
        \includegraphics[width=\linewidth]{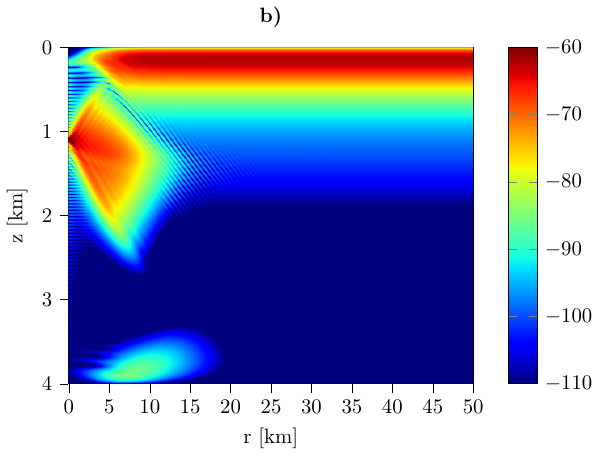}
    \end{subfigure}
        \vspace{0.2cm}

    \begin{subfigure}{0.48\textwidth}
        \centering
        \includegraphics[width=\linewidth]{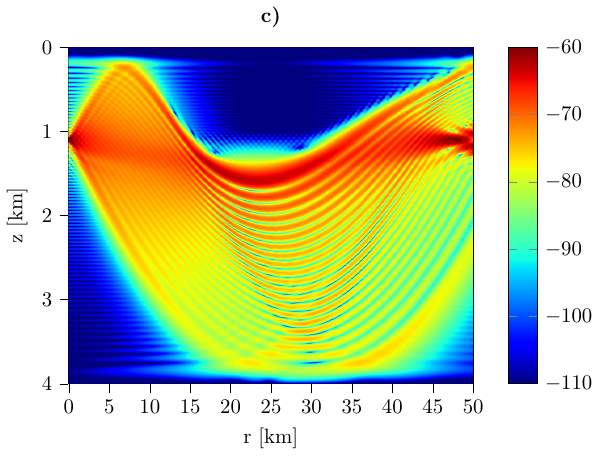}
    \end{subfigure}
    \hfill
    \begin{subfigure}{0.48\textwidth}
        \centering
        \includegraphics[width=\linewidth]{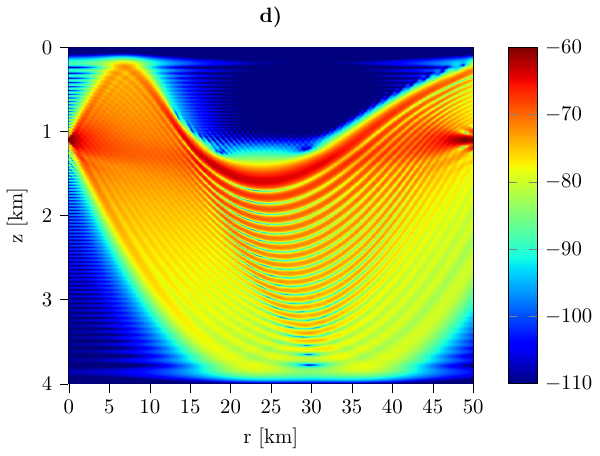}
    \end{subfigure}
    \caption{Acoustic pressure (in dB re $1$~m) computed at $100~\mathrm{Hz}$ using the $[4/4]$ Pad\'e approximation and a step size of $100~\mathrm{m}$ for a deep-water sound channel with the Munk profile~\eqref{eq:Munk}. 
    The subplots (\textbf{a-c}) illustrate the convergence of the Neumann Series~\eqref{eq:neumannseries}.
    \textbf{a)} corresponds to $M=0$, equivalent to the isovelocity case. \textbf{b)} shows the field for $M=1$, \textbf{c)} for $M=2$. \textbf{d)} shows a reference field, generated by the method of normal modes while a slice (the acoustic pressure at a constant depth $z=1300$~m) is displayed in Figure \ref{fig:slice}.
    }
    \label{fig:order}
\end{figure}

\begin{figure}
    \centering
    
    \begin{subfigure}{0.7\textwidth}
        \centering
        \includegraphics[width=\linewidth]{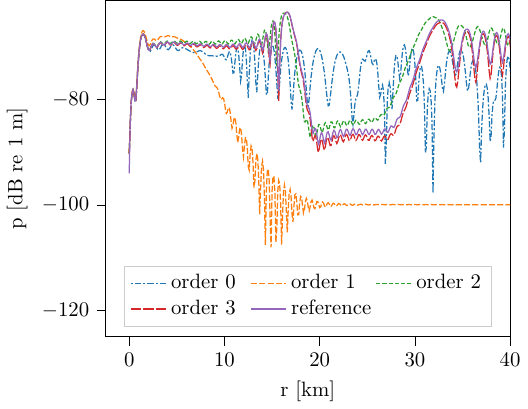}
    \end{subfigure}

    \caption{Acoustic field (in dB re 1~m) propagated through a range independent Munk profile using different correctors of order $M$ at the depth $z=1300~\mathrm{m}$. 
    The reference is once again produced by the method of normal modes. 
    The stepsize used is $h=100~\mathrm{m}$, the discretization constant in depth is $\Delta z = 0.5~\mathrm{m}$.}
    \label{fig:slice}
\end{figure}

Both the numerical experiments and the reference solution depicted in Figure~\ref{fig:order} use 8001 equidistant points in depth and a step size of $h=100$~m.
As can be seen in the comparison, including only three terms of the series ($M=2$) yields results that are already qualitatively comparable to the normal modes reference solution.
However, a difference between the SSSP solution with $M=2$ and the analytical solution is still noticeable, especially at large depths where the perturbation $\delta k^2$ becomes large.
When only the first two terms of the series are included, as in Figure~\ref{fig:order} \textbf{b)}, the computed field exhibits an unphysical blowup due to the first term over-correcting in the positive direction. 
This has been adjusted by renormalizing the field to properly display its shape.
Figure~\ref{fig:slice} displays a slice of the same field at a constant depth $z=1300$~m to show the convergence of this method.
The order 1 field displayed here is normalized once again.
One can also see the alternating sign of the Neumann series here, as the method converges in an alternating fashion from below and above towards the analytical solution.

\subsection{Range Dependent Propagation}
To demonstrate the method's range-dependent capabilities, we introduce a perturbation of the sound speed profile in the form of a synoptic eddy with an exaggerated magnitude. The eddy is parameterized using the formula
\begin{align}
    c_{\rm Eddy}(r,z) &= - c_m \exp\biggl( - \frac{(x - x_0)^2}{r^2_x} \biggr) \exp\biggl( - \frac{(y - y_0)^2}{r^2_y} \biggr) \\
    &\times \biggl( \frac{z - z_0}{r_z} \biggr) \exp\biggl( - \frac{\beta(z - z_0)^2}{r^2_z} \biggr) \,,
    \label{eq:Eddy}
\end{align}
obtained by fitting the CTD measurements of a real eddy in the Sea of Japan \cite{sorokin2024parameterization}.
The parameters used to yield a noticeable perturbation are given in Table~\ref{tab:eddy}.
Combining~\eqref{eq:Munk} and~\eqref{eq:Eddy}, we express the total sound speed distribution by
\begin{equation}
    c(r, z) = c_{\rm Munk}(z) + c_{\rm Eddy}(r, z) \,.
\end{equation}
In ocean acoustics, warm synoptic eddies defocus the field.
While this effect is measurable in a natural eddy, it is not as noticeable in two-dimensional image plots.
For this reason, the perturbation used in this work is exaggerated by a factor of slightly more than ten.
The acoustic field generated by a point source located at $z = z_s = 1100~\mathrm{m}$ is computed using the classical finite-difference SSP method and the SSSP method proposed here. 
The two solutions are then compared to each other and to the acoustic field in the unperturbed sound channel (i.e., without the eddy).

To compute higher-quality reference and SSSP solutions, a higher order ($[6,6]$) Padé approximation is used.
The depth dimension is discretized into 2048 equidistant points, and the step size in the range is $h=100$~m.

\begin{table}[]
    \centering
    \begin{tabular}{cccccccc} \hline
        $\beta$ & $r_x$ & $r_y$ & $r_z$ & $x_0$ & $y_0$ & $z_0$ & $c_m$  \\ \hline
         $1.7125$ & $32~\mathrm{km}$ & $18~\mathrm{km}$ & $250~\mathrm{m}$ & $50~\mathrm{km}$ & $0~\mathrm{m}$ & $1100~\mathrm{m}$ & $40~\mathrm{ms^{-1}}$ \\ \hline
    \end{tabular}
    \caption{Eddy parameters based on the parametrization given in \cite{sorokin2024parameterization}.}
    \label{tab:eddy}
\end{table}

\begin{figure}[ht!]
    \centering
    \begin{subfigure}{0.48\textwidth}
        \centering
        \includegraphics[width=\linewidth]{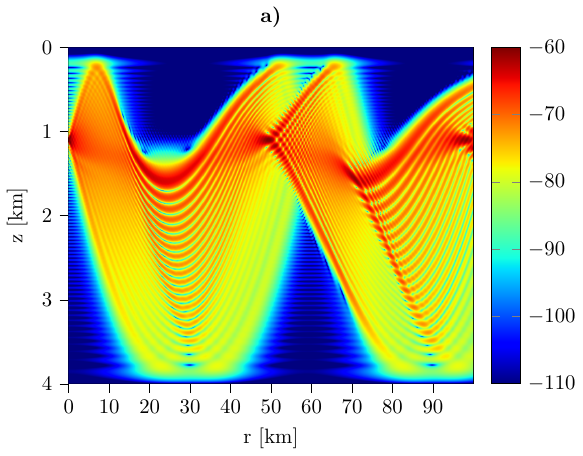}
    \end{subfigure}
    \hfill
    \begin{subfigure}{0.48\textwidth}
        \centering
        \includegraphics[width=\linewidth]{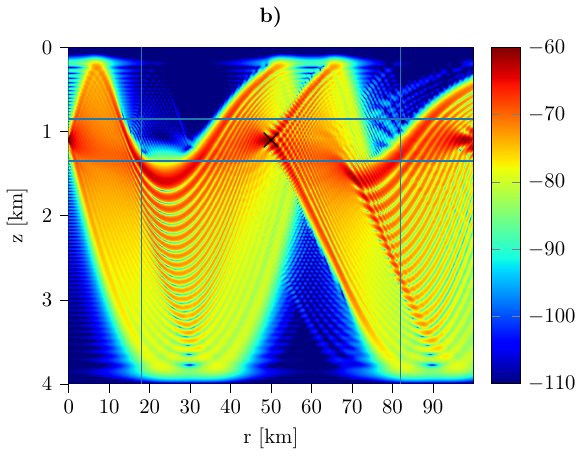}
    \end{subfigure}
    \vspace{0.5cm}
    \begin{subfigure}{0.48\textwidth}
        \centering
        \includegraphics[width=\linewidth]{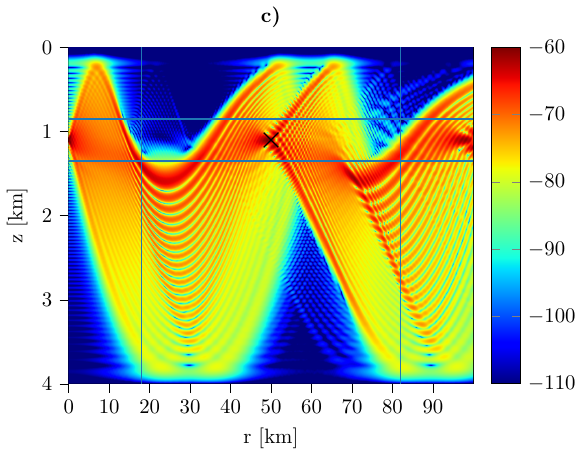}
    \end{subfigure}
    \hfill
    \begin{subfigure}{0.48\textwidth}
        \centering
        
        \includegraphics[width=\linewidth]{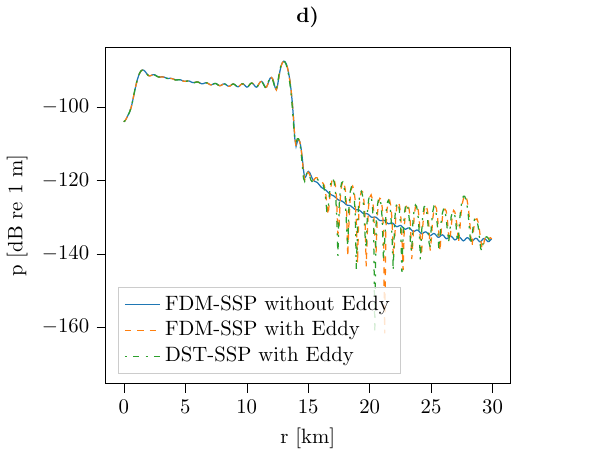}
    \end{subfigure}
    \caption{Acoustic pressure (in dB re 1~m) due to a point source in the problem of sound propagation through a synoptic eddy. 
    In all SSP-simulations, a $[6,6]$ Pad\'e approximation, a step size of $h=100~\mathrm{m}$ and 2048 equidistant discretization points in depth are used. In \textbf{a)} the field without the eddy is displayed as a reference, computed using SSP. 
    The field perturbed by a synoptic eddy computed using SSP is shown in \textbf{b)}.
    The eddy center is indicated by a black ``X'', and its sphere of influence, measured by its variance, is indicated by the horizontal and vertical blue lines. 
    \textbf{c)} shows the acoustic field propagated by SSSP through the eddy.
    \textbf{d)} shows the acoustic field at $z=900\mathrm{m}$, as the field enters the perturbation in sound speed caused by the eddy.}
    \label{fig:eddy} 
\end{figure}
It is clear that the SSSP method provides an approximation of the acoustic field that is just as accurate as the one obtained by the classical SSP method.
A slight difference between the two methods appears in the sliced field at $z=900$~m, but the maximum occurs at the minimum of the acoustic pressure, where accuracy is less important.

\subsection{Coarse Grid Propagation}
To demonstrate numerically that the SSSP algorithm employs an analytical representation of the derivative operator, 
we propagate coarse-grid modal starters using both SSSP and SSP.
Typically, the SSP algorithm begins to exhibit unphysical behavior as soon as the starter cannot be smoothly represented on the discrete grid.
While this is not usually a problem for most starters, which are typically sufficiently smooth, Greene's starter~\eqref{eq:Greene} features a sharp, localized peak about the size of a single wavelength. 
Thus, it requires a large number of points, typically more than 10, per wavelength.
For spectral methods, however, this lower limit is much less strict.
Aliasing effects, the main source of error in discrete transform-based methods, only appear when the highest-frequency oscillations cannot be represented.
This occurs at approximately two points per wavelength.

To demonstrate the magnitude of this effect, we compare the propagation of the modal starter for the Munk profile with 128, 256, and 512 discretization points between SSSP and SSP.
This corresponds to discretization constants of $h=31.25~\mathrm{m}$ and $h=15.625~\mathrm{m}$, respectively.
With this level of discretization, Greene's starter can no longer be properly represented, and even with the smooth modal starter, FDM exhibits significant errors.
This discretization is not fine enough to resolve the highest-frequency features of the field; 
however, the SSSP method preserves the overall shape of the solution.
\begin{figure}[ht!]
    \centering
    \begin{subfigure}{0.48\textwidth}
        \centering
        \includegraphics[width=\linewidth]{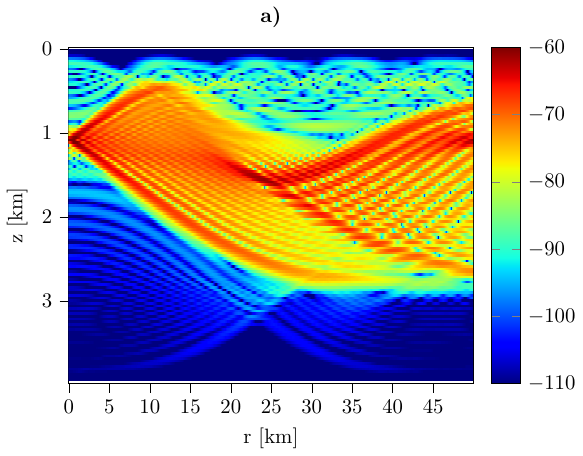}
    \end{subfigure}
    \hfill
    \begin{subfigure}{0.48\textwidth}
        \centering
        \includegraphics[width=\linewidth]{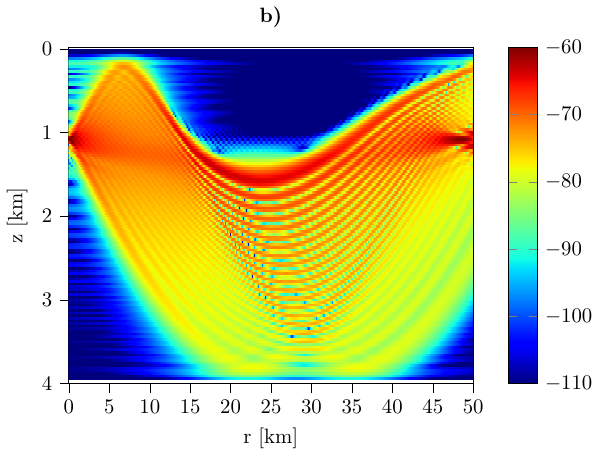}
    \end{subfigure}
    
    \vspace{0.5cm}
    
    \begin{subfigure}{0.48\textwidth}
        \centering
        \includegraphics[width=\linewidth]{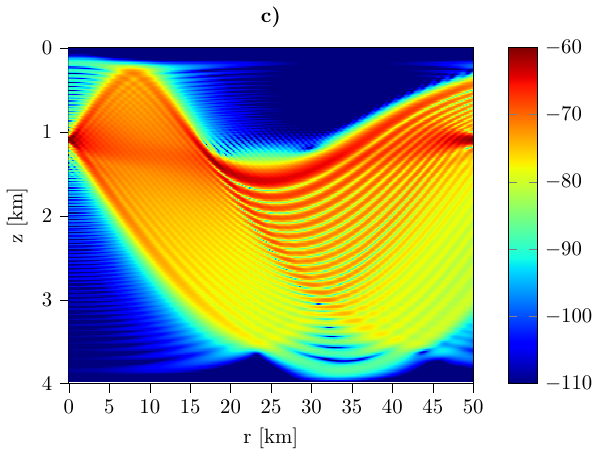}
    \end{subfigure}
    \hfill
    \begin{subfigure}{0.48\textwidth}
        \centering
        \includegraphics[width=\linewidth]{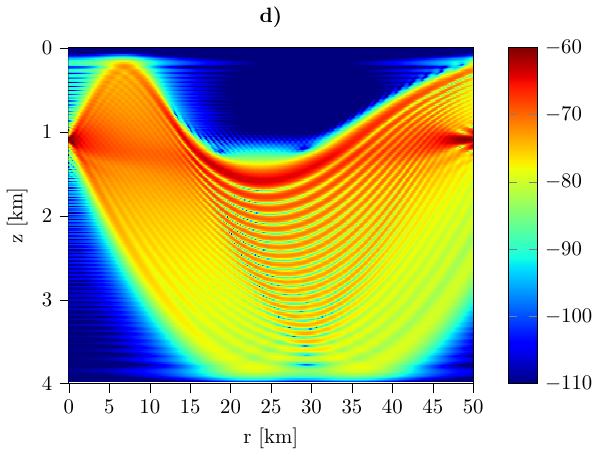}
    \end{subfigure}

    \vspace{0.5cm}
    
    \begin{subfigure}{0.48\textwidth}
        \centering
        \includegraphics[width=\linewidth]{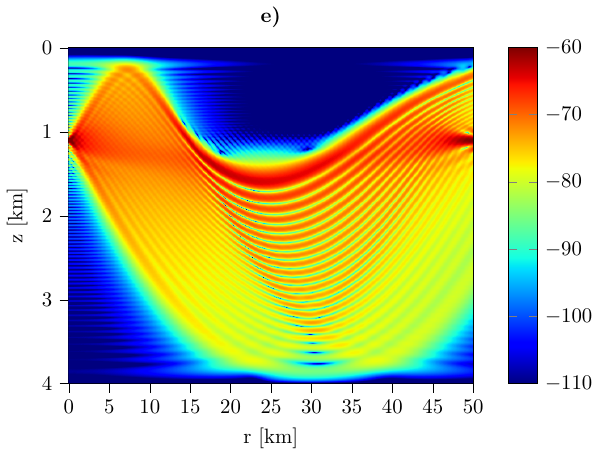}
    \end{subfigure}
    \hfill
    \begin{subfigure}{0.48\textwidth}
        \centering
        \includegraphics[width=\linewidth]{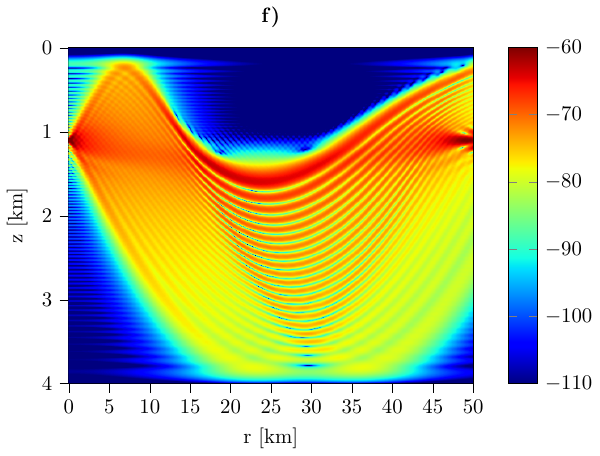}
    \end{subfigure}
    \caption{Propagation of a modal starter through a Munk profile using different discretizations in $z$. 
    Displayed is the acoustic pressure (in dB re 1~m). \textbf{a)} uses SSP with 128 points, \textbf{b)} uses SSSP with again 128 points.
    \textbf{c)} uses 256 points and SSP, while \textbf{d)} uses SSSP. \textbf{e)} uses 512 points and SSP, \textbf{f)} SSSP. 
    A reference solution computed using normal modes is given in Figure~\ref{fig:order}.}
    \label{fig:discretization} 
\end{figure}
Figure~\ref{fig:discretization} shows that SSSP requires only 128 points to converge to the solution for the chosen example.
SSP, on the other hand, has not quite converged, even when using 512 points -- four times as many as SSSP requires.

\section{Conclusion and Outlook}\label{SC}
In this work, we presented a spectral split-step Padé (SSSP) method for solving the one-way Helmholtz equation in wave propagation modeling problems.
This method employs an exact representation of the depth-dependent differential operator using Fourier transforms.
This novel approach eliminates the need for a dense grid in finite difference discretization because no discretization error is introduced by approximating the differential operators.
The SSSP method accommodates standard boundary conditions by employing discrete trigonometric transforms, such as the discrete sine transform.
As with all Fourier transform-based methods, applying the method to artificially truncated computational domains using perfectly matched layers (PMLs) or absorbing boundary conditions, as well as to domains with mixed boundary conditions, is currently not straightforward.

For future work, we will consider combining the pseudospectral approach of Antoine and coauthors \cite{antoine2021derivation, antoine2019simple, antoine2022pseudospectral} with our Fourier method and the use of PML to accomplish the truncation of the computational domain.

\bibliographystyle{plain}
\bibliography{references}

\end{document}